\documentclass[final]{amsart}
\providecommand{\revauthor}[1]{}

\usepackage[notref,notcite]{showkeys}

\usepackage{xspace}
\usepackage{ifthen}
\newcommand{\Ifempty}[3]{\ifthenelse{\equal{#1}{}}{#2}{#3}}

\usepackage{enumerate}
\usepackage[bookmarks,breaklinks,colorlinks,unicode]{hyperref}

\newcommand{\newthm}[2]{\newtheorem{#1}[lemma]{#2}}
\theoremstyle{plain}
\newtheorem{lemma}{Lemma}
\newthm{Ass}{Assumption}
\newthm{convention}{Convention}
\newthm{criterion}{Criterion}
\newthm{fact}{Fact}
\newthm{cor}{Corollary}
\newthm{claim}{Claim}
\newthm{thm}{Theorem}
\newthm{prop}{Proposition}

\theoremstyle{definition}
\newthm{definition}{Definition}
\newthm{example}{Example}

\theoremstyle{remark}
\newthm{remark}{Remark}
\newthm{notation}{Notation}

\newcommand{\Label}[2]{\label{#1:#2}}

\newcommand{\Emph}[1]{\emph{#1}\index{#1}}
\newcommand{\ra}{\rightarrow}

\newcommand{\DefAlias}[2]{\expandafter\xdef\csname #1\endcsname{#2}}
\newcommand{\CiteAlias}[2]{\DefAlias{CITE#1}{#2}}
\newcommand{\Cite}[1]{\cite{\csname CITE#1\endcsname}}

\CiteAlias{ACFA}{MR1652269}
\CiteAlias{LAG}{MR1102012}
\CiteAlias{DG}{MR1972449}
\CiteAlias{groupoid}{groupoid}
\CiteAlias{bounds}{MR739626}
\CiteAlias{stability}{MR1416106}
\CiteAlias{SGA}{MR0354652}
\CiteAlias{singer}{MR1480919}
\CiteAlias{milne}{MR559531}
\CiteAlias{topoi}{MR1300636}
\CiteAlias{sacks}{MR0398817}
\CiteAlias{EI}{MR1083551}
\CiteAlias{HrWeil}{MR1827833}
\CiteAlias{SingerMod}{SingerMod}
\CiteAlias{PillayZoe}{MR1650667}
\CiteAlias{Bouscaren}{MR984628}
\CiteAlias{eisenbud}{MR1322960}

\newcommand{\Email}{\url{mailto:m.kamensky@uea.ac.uk}}

\newcommand{\x}{\times}
\newcommand{\MapsTo}{\rightarrow}
\newcommand{\df}{\stackrel{\mathrm{def}}{=}}
\newcommand{\w}{\mbox{${\omega}$}}
\newcommand{\ti}[1]{\tilde{#1}}
\newcommand{\bigand}[1]{{\bigwedge}_{#1}}
\newcommand{\Pspace}{\mathbb{P}}
\newcommand{\ACF}{$ACF$}

\newcommand{\Th}{\mbox{$\mathfrak T$}}
\def\.{\cdot}
\def\:{\colon}
\newcommand{\E}{\exists}
\newcommand{\tensor}{\otimes}
\newcommand{\emb}{\hookrightarrow}
\newcommand{\p}{\mathfrak{p}}
\newcommand{\q}{\mathfrak{q}}
\def\A{\mbox{$\mathbb A$}}
\newcommand{\LL}{\mathcal{L}}
\newcommand{\Z}{\mathbb{Z}}
\newcommand{\Gen}[1]{{\langle #1 \rangle}}
\newenvironment{Thanks}{\subsection*{Acknowledgement}}{}

\title[model completion of modules]{The model completion of the theory of 
modules over finitely generated commutative algebras}
\author{Moshe Kamensky}
\revauthor{Kamensky, Moshe}
\address{Department of Mathematics\\
         The Hebrew University\\
         Jerusalem, Israel}
\curraddr{Department of Maths\\
         University of East Anglia\\
         Norwich, NR4 7TJ, England}

\email{\Email}
\urladdr{\url{http://mkamensky.notlong.com}}


\begin{document}


\begin{abstract}
  We find the model completion of the theory modules over $\mathbb{A}$, where 
  $\mathbb{A}$ is a finitely generated commutative algebra over a field $K$.  
  This is done in a context where the field $K$ and the module are 
  represented by sorts in the theory, so that constructible sets associated 
  with a module can be interpreted in this language. The language is expanded 
  by additional sorts for the Grassmanians of all powers of $K^n$, which are 
  necessary to achieve quantifier elimination.

  The result turns out to be that the model completion is the theory of a 
  certain class of ``big'' injective modules. In particular, it is shown that 
  the class of injective modules is itself elementary. We also obtain an 
  explicit description of the types in this theory.
\end{abstract}

\maketitle

\section{Introduction}

An algebra $\A$ finitely generated over an algebraically closed field $K$ 
corresponds to an affine variety $V$ over $K$, and a module over $\A$ 
corresponds to a (quasi-coherent) sheaf over $V$. Whereas varieties can be 
reasonably considered within the framework of model theory (for example, as 
definable sets in the theory \ACF{} of algebraically closed fields), modules 
(or sheaves) do not appear so naturally. For example, basic results about 
definability of the fibre dimension are proved, using algebraic methods, for 
algebras and modules alike. On the model theoretic side, the fibre dimension 
for a map of varieties (or more generally, for definable sets) is well 
understood, in a much more general framework. However, the analogous 
statements for modules can not even be phrased. This work represents, we 
hope, a first step in approaching these questions.

The purpose of this paper is to find the model completion for the theory of 
modules over a finitely generated commutative $K$-algebra ($K$ a field), and 
describe the types in that theory. Our initial approach in formulating this 
theory is to use a two-sorted language, with a sort $K$ for the field, and 
another sort $M$ for the module. In addition to the field structure on $K$ 
and the $K$-vector space structure on $M$, we introduce symbols for $n$ 
commuting linear operators on $M$, that represent the generators of the 
algebra.

Our goal is to find the model completion. To estimate the feasibility of our 
task, we consider the case $n=0$. In this case we simply have a vector space 
$M$ over $K$. We immediately observe that the most basic relation on $M$, 
that of linear dependence, can not be expressed in this theory without 
quantifiers.  This example leads us to introduce additional sorts for all the 
Grassmanians\index{Grassmanian} of the vector spaces $K^n$. The dependence 
relation on $M$ then takes values in these Grassmanians. Thus, with this 
addition to the language, the above problem is resolved, and it turns out 
that this is the only obstacle for the existence of a model completion, even 
for the case when $n>0$.

Below we give a precise definition of the language and the theory we work 
with. The rest is divided into the cases $n=0$ and $n>0$. Although the first 
case is not really different, the kind of problems the cases deal with are 
different and independent: in the first case, we deal with the vector space 
structure, as well as the new sorts introduced. In the second case, the main 
interest comes from the action of the operators $T_i$.

\begin{Thanks}
This work is part of my PhD research, performed in the Hebrew university 
under the supervision of Ehud Hrushovski. I would like to thank him for his 
guidance, and in particular for suggesting this question, and helping with 
the difficulties, as they arose.
\end{Thanks}

\subsection{The theory of modules over a commutative $K$-algebra}

Given $n \geq 0$, we use the following language 
$\LL=\LL_n$:\index{aa@$P_{\ti{\varphi}}$|(}
\begin{align}
  \LL_n = (& K, +, \., 0, 1, \label{eqn:Lfield}\\
  & M, +, 0, \., \label{eqn:Lmodule}\\
  & G^i_{(i > 0)}, {\pi_i}_{(i > 0)}, P_{\tilde{\varphi}}, {D_i}_{(i>0)} 
  \label{eqn:Lgrass}\\
  & T_1,\dots,T_n )\label{eqn:Loperators}
\end{align}
Where
\begin{itemize}
  \item (\ref{eqn:Lfield}) is the language of fields
  \item (\ref{eqn:Lmodule}) is the language of abelian groups together with a 
    function symbol $\.\:K\x M\MapsTo M$
  \item each $G^i$ is a sort, and $\pi_i\:K^{i^2}\MapsTo G^i$ is a function 
    symbol
  \item $\tilde{\varphi}$ is a quantifier free formula 
    $\varphi(x_1,\dots,x_N)$ in the language of fields, together with a 
    partition of its variables to sets of sizes $k_1^2,\dots,k_m^2, l$.  
    Given such a $\tilde{\varphi}$, $P_{\tilde{\varphi}}$ is a predicate 
    symbol on $G^{k_1}\x\dots\x G^{k_m}\x K^l$ (For the meaning of all this, 
    see below)
  \item Each $D_i:M^i\MapsTo G^i$ is a function symbol
  \item Each $T_i$ is a function symbol $T_i:M\MapsTo M$
\end{itemize}

Given an ideal $I\subset \Z[T_1,\dots,T_n]$, the theory 
$\Th=\Th_I$\index{aa@$\Th$}\index{aa@$\Th_I$} says the following:
\begin{itemize}
  \item
    (\ref{eqn:Lfield}) is a field, and (\ref{eqn:Lmodule}) is a vector space 
    over it

  \item
    For each $i$, $(G^i, \pi_i)$ is the set of all linear subspaces (the 
    \Emph{Grassmanian}) of $K^i$.  Thus, we view $K^{i^2}$ as an $i$-tuple of 
    row vectors in $K^i$, and the theory says that $\pi_i:K^{i^2}\ra G^i$ is 
    surjective, and two elements belong to the same fibre if and only if the 
    corresponding row vectors span the same linear subspace.

  \item
    For any $\tilde{\varphi}$, $P_{\tilde{\varphi}}$ is (using the notation 
    above) the set induced on $G^{k_1}\x\dots\x G^{k_m}\x K^l$ by $\varphi$, 
    i.e., we have the formula
    \begin{equation*}
      \begin{split}
      & \forall p_1\in G^{k_1},\dots,p_m\in G^{k_m},\bar{x}\in K^l \\
      & \phantom{\forall p_1\in}(P_{\tilde{\varphi}}(p_1,\dots,p_m,\bar{x}) 
      \iff\\
      & \phantom{\forall p_1\in(} \E \bar{x_1}\in 
      K^{k_1^2},\dots,\bar{x_m}\in K^{k_m^2} 
      (\varphi(\bar{x_1},\dots,\bar{x_m},\bar{x}) \land \\
      &  \phantom{\forall p_1\in( \E \bar{x_1}\in 
      K^{k_1^2},\dots,\bar{x_m}\in K^{k_m^2}(} 
      \bigwedge_i\pi(\bar{x_i})=p_i))
    \end{split}
    \end{equation*}
    Note that this is not an additional structure, but part of $K^{eq}$. 
    However, as explained above, this addition (together with the operators 
    $D_i$) is essential to achieve quantifier 
    elimination.\index{aa@$P_{\ti{\varphi}}$|)}

  \item
    For any $\bar{v}=v_1,\dots,v_m\in M$, $D_m(\bar{v})$ is the subspace 
    $p\in G^m$ of all $\bar{x}\in K^m$ such that $\sum x_iv_i=0$.

  \item
    The $T_i$ represent generators of the algebra. Thus they are commuting 
    linear operators on $M$, and any $p(T_1,\dots,T_n)\in I$ is the $0$ 
    operator (See remark~\ref{rmk:noncom} about the non-commutative case.)

\end{itemize}

Models of this theory are (determined by) pairs $(K,M)$, where $K$ is a 
field, and $M$ is a module over $K[T_1,\dots,T_n]/{\hat{I}}$. Here, 
$\hat{I}$ is the ideal generated by $I$ in $K[T_1,\dots,T_n]$.

\begin{notation}
  For the sake of readability, we use the following conventions: the letter 
  $G$ is used to denote a $G^i$ with an unspecified $i$. Unless otherwise 
  mentioned, $x,y,z$ are field variables, $p,q,r$ are $G$ variables, and 
  $u,v,w$ are module variables. $X,Y,Z$ are used for tuples of field 
  variables when considered as matrices. Also, since the number of operators 
  $T_i$ is fixed in every situation, $n$ is released for other uses.
\end{notation}

\section{The case $n=0$}\Label{sec}{eq0}

We are looking for a model completion of the theory above, so in particular 
the field part $K$ should eliminate quantifiers. Since in $\LL$, the 
quantifier free subsets of the field are only those defined in the field 
language, this leads us to the requirement that $K$ is algebraically closed.  
This requirement makes the theory complete, up to the characteristic of $K$ 
and the dimension of $M$ as a vector space over $K$ (in fact, fixing the 
characteristic of $K$ and the dimension of $M$, the theory we get is 
$\aleph_1$ categorical). Let $\tilde{\Th}$\index{aa@$\ti{\Th}$} be $\Th_0$, 
together with the axioms saying that $K$ is algebraically closed, and $M$ is 
of a given dimension over $K$ (which might be infinity). Our first goal is:

\begin{prop} \label{prp:eq0}
  The theory $\tilde{\Th}$ has elimination of quantifiers
\end{prop}

We begin with a few remarks concerning only the relation between $K$ and the 
$G^i$. For $\varphi(\bar{x},p_1,\dots,p_k)$ a formula, let

\begin{equation*}
  \varphi^*(\bar{x},Y_1,\dots,Y_k)\df 
  \varphi(\bar{x},\pi(Y_1),\dots,\pi(Y_k))
\end{equation*}

Such formulas will be called homogeneous (in $Y_1,\dots,Y_k$).

For $p_i\in G^{l_i}$ ($1\le i\le k$), let $[p_1,\dots,p_k]\in 
G^{\sum l_i}$ be the subspace $\oplus p_i$ of $\oplus K^{l_i}$.

We are going to make assertions regarding linear transformations on 
spaces like $K^n$ and $M^n$. Since we usually need the claims for 
matrices with arbitrary terms (and not just constants), we redefine 
`linear map' to mean any definable map from $K^m$ to the set of 
$n_1\times n_2$ matrices (some $m,n_i$), considered as linear 
transformations acting on the right for $K$ and on the left for $M$.

\begin{lemma} \label{lma:grass}
  The following facts hold in $\tilde{\Th}$:
  \begin{enumerate}[a.]
    \item Any quantifier free formula is equivalent to a quantifier 
      free formula without $\pi$.\label{lma:grass1}
    \item Let $\varphi(\bar{x},p_1,\dots,p_k)$ be a quantifier free formula. 
      Then
      \begin{equation} \label{eqn:star}
        \varphi\equiv P_{\varphi^*}
      \end{equation} \label{lma:grass2}
    \item The map $(p_1,\dots,p_k)\mapsto [p_1,\dots,p_k]$ is definable 
      without quantifiers.
    \item The theory $\tilde{\Th}$ restricted to $K$ and the $G^i$ eliminates 
      quantifiers. Thus, the equation~\eqref{eqn:star} holds for any formula 
      $\varphi$ in the restricted theory. Any formula in this theory is 
      equivalent to the formula $P_\varphi$ for some field formula $\varphi$.
    \item For any linear map $A:K^m\MapsTo K^n$, and for subspaces 
      $p\subseteq K^m$ and $q\subseteq K^n$, the image $pA$ and the 
      inverse image $qA^{-1}$ are again linear subspaces.  Thus we 
      have induced maps $A:G^m\MapsTo G^n$ and $A^{-1}:G^n\MapsTo G^m$ (these 
      maps will be written on the left).
  \end{enumerate}
\end{lemma}

\begin{proof}
  \begin{enumerate}[a.]
    \item It's enough to prove this for atomic formulas.  There are 
      two kinds of these:
      \begin{itemize}
        \item $P_{\varphi}(\pi(\bar{x}),\dots)$.  This holds, by definition, 
          iff
          \begin{equation*}
            \E\bar{y},\dots(\varphi(\bar{y},\dots)\land 
            \pi(\bar{y})=\pi(\bar{x})\land\dots)
          \end{equation*}
          where $\bar{y}$ does not appear in any of the $\dots$ parts. The 
          expression $\pi(\bar{y})=\pi(\bar{x})$ just says that $\bar{x}$ and 
          $\bar{y}$ span the same vector space, which can be expressed using 
          only the language of fields.  Therefore, the above formula is 
          equivalent to $\E\dots(\E \bar{y}(\varphi'(\bar{x},
          \bar{y}, \dots))\land\dots)$
          where $\varphi'$ is in the language of fields. By elimination of 
          quantifiers in \ACF{}, $\E \bar{y}(\varphi'(\bar{x},\bar{y}, 
          \dots))$ is equivalent to some quantifier free $\varphi''$, so our 
          original formula is equivalent to $P_{\varphi''}$, where we got rid 
          of one $\pi$.
        \item $\pi(\bar{x})=p$. This one is equivalent to 
          $P_{\bar{x}=\bar{y}}(\bar{x},p)$, with the corresponding partition 
          of the variables.
      \end{itemize}

    \item For $\varphi=P_\psi$, this follows directly from the definitions.  
      By \ref{lma:grass1}, this is the only kind of atomic formulas we should 
      check.  On the other hand, $*$ is a homomorphism of boolean algebras, 
      and so is $P$ restricted to formulas of the form $\psi^*$, so the 
      result follows.

    \item This map is $P_{\bar{x}=\bar{y}}$ for an appropriate partition of 
      the variables, and with $\bar{x}$ padded with zeroes in the right 
      places (more precisely, $\bar{x}$ is a matrix with $k$ matrices of the 
      right sizes on the diagonal, and $0$ elsewhere).

    \item By~\ref{lma:grass2}, we need to show there is a quantifier free 
      formula equivalent to
      \begin{equation*}
        \E A_1(P_\varphi(\bar{x},p_1,\dots,p_k))
      \end{equation*}
      where $A$ is either $x$ or $p$. Unravelling the definition, this 
      amounts to the fact that existential quantifiers commute, together with 
      quantifier elimination for algebraically closed fields. The rest is 
      just a summary of the previous items, together with the fact (clear 
      from inspecting the proofs) that they can be put together.

    \item Is obvious.
  \end{enumerate}

\end{proof}

For $A_1,\dots,A_k$ linear maps, and $\varphi(\bar{x},p_1,\dots,p_k)$ a 
formula, we set\index{aa@$A_*$}\index{aa@$A^*$}
\begin{equation*}
(A_1,\dots,A_k)^*\varphi=\bar{A}^*\varphi\df \varphi(\bar{x}, 
A_1(q_1), \dots, A_k(q_k))
\end{equation*}
and
\begin{equation*}
  (A_1,\dots,A_k)_*\varphi=\bar{A}_*\varphi\df \varphi(\bar{x}, 
  A_1^{-1}(q_1), \dots, A_k^{-1}(q_k))
\end{equation*}

Note that these formulas will depend on the additional variables of 
the $A_i$.\footnote{
  Strictly speaking, the operators $A_i$ do not actually exist in the 
  language, but the formulas exist (and are quantifier free).
}

We now go back to the full $\tilde{\Th}$, and the next step is to 
analyse the quantifier free formulas in the theory. The main lemma we 
need is:

\begin{lemma} \label{lma:main}
  Any quantifier free formula $\varphi$ is equivalent to a quantifier 
  free formula $\psi$, in which for every term of the form 
  $D(t_1,\dots,t_n)$, each $t_i$ is either a module variable or a 
  module constant\footnote{
    We assume that the base set is a substructure
  }
\end{lemma}

\begin{proof}
  The claim follows from the fact that for any linear map $A$, we have 
  $D(A\bar{v})=A^{-1}(D(\bar{v}))$: Indeed, both sides are equal to the space 
  of all $\bar{x}$ such that $\bar{x}A\bar{v}=0$. Now, any quantifier free 
  formula $\varphi$ has the form
  \begin{equation*}
    \varphi'(\bar{x}, D(A_1\bar{v}),\dots, D(A_k\bar{v}),\bar{p})
  \end{equation*}
  so by the above equality, $\varphi$ is equivalent to
  \begin{equation*}
    (\bar{A}_*\varphi')(\bar{x},D(\bar{v}),\dots,D(\bar{v}),\bar{p})
  \end{equation*}
\end{proof}

We can now prove quantifier elimination:
\begin{proof}[Proof of proposition \ref{prp:eq0}]
  Let $\varphi(\bar{x}, \bar{p}, \bar{v})$ be some quantifier free 
  formula.  We need to find a quantifier free formula equivalent to 
  $\E A\varphi$, where $A$ is one of $x_0$, $p_0$ or $v_0$.  
  Now, by lemma \ref{lma:main}, there is some formula 
  $\varphi'(\bar{x}, \bar{p}, q)$ such that $\varphi$ is equivalent 
  to $\varphi'(\bar{x}, \bar{p}, D(\bar{v}))$.  Hence, for the cases 
  that $A$ is either $x_0$ or $p_0$, $\E A\varphi$ is equivalent 
  to
  $\E A \varphi'(\bar{x}, \bar{p}, D(\bar{v}))$, 
  and $\E A\varphi'$ is equivalent to a quantifier free formula 
  by lemma \ref{lma:grass}. Thus the only case left is 
  $\E v_0\varphi$.

  Let $\varphi'_0$ be the formula
  \begin{equation*}
    \varphi'(\bar{x},\bar{p},q)\wedge
    \E\bar{y}( (1, \bar{y}) \in q)
  \end{equation*}
  (i.e., the set of $q$ satisfying $\varphi$ whose projection to the first 
  coordinate is not $0$), and let
  $\varphi'_1=\varphi'\wedge\neg\varphi'_0$. Since existential quantifiers 
  commute with disjunction, it is enough to prove for each of these cases 
  separately.

  Assume first that $\varphi'=\varphi'_0$. Then
  $\E v_0\varphi$ is equivalent to
  \begin{equation*}
    \E\bar{y}
    \varphi'(\bar{x},\bar{p},D(-\sum_{i>0}{y_iv_i},v_1,\dots,v_n))
  \end{equation*}
  
  In the case $\varphi=\varphi'_1$, $\varphi$ says that $v_0$ is independent 
  of the other vectors, and therefore $D(\bar{v})$ coincides with 
  $i(D(v_1,\dots,v_n))$, where $i:K^n\mapsto K^{n+1}$ is the inclusion as the 
  last $n$ coordinates.  Hence
  $\E v_0\varphi$ is equivalent to
  \begin{equation*}
    \varphi'(\bar{x},\bar{p},i(D(v_1,\dots,v_n)))\wedge
    (<v_1,\dots,v_n>\ne M)
  \end{equation*}
  Since the theory determines the dimension of $M$, the statement that the 
  $v_i$ span $M$ depends only on the dimensions of $D$ applied to subsets of 
  the $v_i$, hence is quantifier free.
\end{proof}

The next goal is to analyse the quantifier free types. Since we don't 
use quantifier elimination here, we will be able to use this to give 
a second proof of quantifier elimination. Then, because of quantifier 
elimination, this will give information about the spaces of types, and 
eventually $\w$-stability will be shown.

To prove quantifier elimination, we will use the following criterion 
(cf~\Cite{sacks}):

\begin{criterion}\Label{crt}{eq}

A theory $T$ eliminates quantifiers if for any model $M$ and any
$A\subseteq M$, any quantifier free $1$-type over $A$ is also a type
(i.e. consistent) with respect to any extension of $T_A$ (where $T_A$ is the 
theory obtained by adding to $T$ all quantifier free sentences over $A$ that 
hold in $M$.)

\end{criterion}

We begin by analysing the substructures of a model of $\tilde{\Th}$, 
and first, as before, we consider only the restriction to the sorts $K$ and 
$G^i$. For this restricted theory we will assume elimination of quantifiers 
(as proved in lemma~\ref{lma:grass}.) Let $A$ be a substructure of a model
of this theory. Then $K(A)$ is an integral domain, whose fraction field we 
denote by $L$, $M(A)$ is a vector space over $L$, and $G^i(A)$ contains all 
the subspaces of $L^i$.

\begin{claim}\Label{clm}{goodG}
  There is a unique minimal extension $B$ of $A$ such that $K(B)$ is 
  a field, and for each $i$, $\pi_i:K(B)^{i^2}\MapsTo G(B)^i$ is onto.
\end{claim}

\begin{proof}
  First, we may assume that $K(A)$ is a field by passing to the fraction 
  field. Consider the subset $\Pspace(A)$ of $G(A)$ consisting of the 
  one-dimensional subspaces. A point of this subset corresponds to a line in 
  some affine space $K^i$. For $\pi_i$ to be onto, this line should have a 
  point in $B$.  This will happen if and only if the unique point on this 
  line whose $k$-th coordinate is $1$ has its other coordinates in $B$.  Such 
  points correspond to intersection of this line with the standard cover of 
  $\Pspace$. This cover corresponds to some elements (over $0$) of $G$, and 
  the intersections are encoded in the structure of $G$. We thus get a finite 
  set of points in affine space, one for each such intersection. Using the 
  standard projections, we get a finite set of points in $\A^1$.  The type of 
  these points as field elements is well defined, since both the field $K(A)$ 
  and the field operations can be viewed as part of the structure $G$. We 
  thus get a field extension $K(B)$, which, by construction, contains a point 
  in each element of $\Pspace(A)$, and is obviously minimal with this 
  property.

  It remains to show that $K(B)$ contains a basis for any other element of 
  $G(A)$ as well. Consider the elements $G_k^n(A)$ in $G(A)$ corresponding to 
  $k$ dimensional subspaces of $K^n$. This set has a natural embedding (over 
  $\Z$) into $\Pspace(\bigwedge^k{}K^n)$, corresponding to the natural map 
  $(K^n)^k\MapsTo\bigwedge^k{}K^n$. For a given point of $G_k^n(A)$, its 
  image in the above projective space contains, by the definition of $B$, a 
  point of $K(B)$. Thus the problem reduces to showing that, for any field 
  $L$, any point of $\bigwedge^k{}L^n$ has a pre-image in $(L^n)^k$ under the 
  natural map.  However, the pre-image set is a $GL_k$-torsor (the action of 
  $GL_k$ corresponds to changing the vector space basis), and any such torsor 
  has an $L$ point (see, e.g., \Cite{milne}, Lemma 4.10.)

\end{proof}

Next, we extend the statement to the sort $M$:

\begin{claim} \label{clm:good}
  Let $A$ be a substructure, and let $B$ be as promised by 
  claim~\ref{clm:goodG}. Then there is a unique minimal vector space $V$ over 
  $K(B)$ such that $(B,V)$ is an extension of $A$ as a substructure.
\end{claim}
\begin{proof}
  Let
  \begin{equation*}
    V=K(B)\tensor_{K(A)}M(A)/
    \Gen{\sum x_i\tensor v_i | \bar{x}\in D(\bar{v})}
  \end{equation*}
  Since
  $D(\sum x^1_i\tensor v^1_i,\dots,\sum x^m_i\tensor v^m_i)$ is 
  determined by $D(\bar{v^1},\dots,\bar{v^m})$, this already defines 
  a structure. It is obvious that this is what we want.
\end{proof}

Let us say that $A$ is a \Emph{good substructure} if $K(A)$ is a 
field and $G^i(A)$ is the set of subspaces of $K(A)^i$ (in other 
words, it is a model of $\Th$). Then the above claims say that any 
substructure has a unique minimal extension to a good substructure (in other 
words, definably closed structures are good.)

\begin{claim} \label{clm:tensor}
  Let $A$ be a good substructure, $B$ an extension of $A$.  Then 
  $M(B)$ contains $K(B)\tensor_{K(A)}M(A)$.
\end{claim}

\begin{proof}
  Since $M(B)$ is a vector space over $K(B)$ containing $M(A)$, there 
  is a canonical map $i:K(B)\tensor_{K(A)}M(A)\MapsTo M(B)$. Assume that 
  $\sum x_j\tensor v_j$ goes to $0$ in $M(B)$ ($v_j\in M(A)$).  Then 
  $\bar{x}\in D(\bar{v})$, but according to the assumption, 
  $D(\bar{v})$ has a basis with coordinates in $K(A)$, so $\sum 
  x_j\tensor v_j$ is $0$ already in $K(B)\tensor_{K(A)}M(A)$.
\end{proof}

This implies that for good substructures, statements regarding the 
vector space are unambiguous: In general, for example, the statement 
``$v_1,\dots,v_n$ are linearly independent'' might mean either that 
it is independent over the field part of the structure, or that 
$D(\bar{v})=0$. For good substructures, this is the same.

We can now give a

\begin{proof} [Second proof of quantifier elimination]
  We use criterion~\ref{crt:eq}. Let $A$ be a substructure. Any model of 
  $\tilde{\Th}$ containing $A$ will also contain the substructure given by 
  claim~\ref{clm:good}, and $K$ will contain an algebraic closure of $K(A)$, 
  so by claim~\ref{clm:tensor} we may assume that
  \begin{itemize}
    \item $K(A)$ is an algebraically closed field.
    \item each $\pi$ is onto.
    \item $M(A)$ is a vector space over $K(A)$.
  \end{itemize}

  Let $\Th_1$ be a theory extending $\tilde{\Th}_A$, and let $\p(v)$ be a 
  quantifier free $1$-type over $A$ with respect to $\tilde{\Th}$, in the 
  module sort (quantifiers on other sorts are eliminated as before). Consider 
  the set of formulas $D(v,v_1,\dots,v_n)\neq 0$ with $v_i\in M(A)$, 
  satisfying $D(v_1,\dots,v_n)=0$. Assume first that there is no such 
  formula in $\p$. Then (since $\p$ is consistent), the vector space is 
  either $\infty$-dimensional, or of dimension greater than the dimension of 
  $M(A)$.  In any case, there is a model of $\Th_1$ which has a member 
  outside the space generated by $M(A)$. Any such member will satisfy the 
  type.

  Now assume, conversely, that there are formulas as above, and assume that 
  $n$ is minimal. Then for any $u_i$ with $D(v,\bar{u})\neq{}0$, the space 
  $V$ spanned by $v_1,\dots,v_n$ is contained in the space spanned by 
  $\bar{u}$ (Otherwise, the intersection of these spaces is properly 
  contained in $V$, and any basis of it is a contradiction to the minimality 
  of $n$). We claim that the set
  \begin{equation*}
    \{P_\varphi(D(v,\bar{v}))\in\p\}
  \end{equation*}
  determines the type. Let
  $\psi(D(v,\bar{u}_1),\dots,D(v,\bar{u}_k))$ be a formula in $\p$. 
  We first note, that if $\bar{w}_1$ spans the same subspace as 
  $\bar{u}_1$, then $\psi$ is equivalent to some formula 
  $\psi'(D(v,\bar{w}_1),\dots,D(v,\bar{u}_k))$: by assumption, there is some 
  matrix $U$ (over $K(A)$!) such that $(v,\bar{u}_1)=U(v,\bar{w}_1)$. Hence
  the equivalence follows from lemma~\ref{lma:main}. In particular, we may 
  assume that the first $n$ vectors in each $\bar{u}_i$ coincide with 
  $\bar{v}$, and that each $\bar{u}_i$ is linearly independent. But then, 
  letting $i_m:K^n\emb K^{l_m}$ be the inclusion of the first $n$ coordinates 
  (where $l_m$ is the length of $\bar{u}_m$), it is clear that 
  $(i_1,\dots,i_k)^*\psi$ is the formula we seek.

  Let $\q=\{P_\varphi(p):P_\varphi(D(v,\bar{v}))\in\p\}$.  Since \ACF{} 
  eliminates quantifiers, there is a model of $\Th_1$ in which $\q$ has a 
  realisation, $q$.  Since  $\bar{v}$ is independent, $q$ will be of 
  dimension either $1$ or $0$.  If it's $1$, let $x_1,\dots,x_n$ be the 
  unique tuple with $(1,\bar{x})\in q$.  Then $v=-\sum x_iv_i$ satisfies 
  $\p$.  If $q=0$, then the dimension of $M$ must be more than $n$ (otherwise 
  $\p$ would be inconsistent, since the dimension is given already in $\Th$).  
  Then any $v$ independent of $\bar{v}$ satisfies $\p$.
\end{proof}

Let's record the result in the proof as a separate claim:
\begin{claim}[description of the vector space types]
  A $1$-type $\p(v)$ over a good structure $A$ is determined by either a 
  sequence $v_1,\dots,v_n\in M(A)$ of minimal length such that 
  \begin{equation*}
    D(v,v_1,\dots,v_n)=0
  \end{equation*}
  is in $\p$, together with the type $\q$ (in the field and Grassmanian 
  sorts) such that
  \begin{equation*}
    \p=\q(D(v, v_1,\dots,v_n))
  \end{equation*}
  or by the fact that there is no such sequence (In other words, it is 
  determined by the minimal subspace to which $v$ belongs, together with the 
  minimal field over which it happens).
\end{claim}

\begin{remark}
  In the proof we dealt only with quantifier free types, but now we 
  know that this is all there is.
\end{remark}

Recall that a theory is \emph{$\w$-stable}\index{stable@$\w$-stable} if the 
set of types over any countable set is countable. As a corollary of the 
description of types we get
\begin{cor}
  $\tilde{\Th}$ is $\w$-stable.
\end{cor}

\begin{proof}
  This follows by counting the types, using the above claim and 
  $\w$-stability of \ACF{}.
\end{proof}

We note that in the case that $M$ is finite-dimensional, this corollary 
already follows from the $\w$-stability of \ACF{}, since, after adding a 
basis, $M$ is interpretable in the field. However, the quantifier elimination 
result holds without adding any parameters.

\section{The general case}

Unlike the case $n=0$, for $n>0$, $\Th_I$ is far from being complete (unless 
$I$ is maximal), even if the field is algebraically closed.  Nevertheless, 
quantifier elimination in the field (and $G$) variables follows automatically 
from the case $n=0$. For the full quantifier elimination, we consider an 
extended theory $\ti{\Th}$\index{aa@$\ti{\Th}$} whose models satisfy the 
following property: Given a model $N$, let $\A$ be the algebra 
$K(N)[T_1,\dots,T_n]/\Gen{I}$.  Then any (finite) set of conditions:
\begin{gather}
  f_iv=v_i\\
  g_jv\notin U_j
\end{gather}
where $f_i$,$g_i$ are in $\A$, $v_i\in{}M(N)$ are module elements, and $U_i$ 
are finite dimensional subspaces of $M(N)$, has a solution $v$, provided 
that:
\begin{itemize}
  \item
    If
    \begin{equation}
      \sum t_if_i=0
    \end{equation}
    then
    \begin{equation}
      \sum t_iv_i=0
    \end{equation}
    for any $t_i\in \A$.
    
  \item
    No $g_i$ is in the ideal generated by the $f_i$.

\end{itemize}

Note, that these conditions are necessary for a solution to exist.

Since, as they are written, these conditions involve quantifying over all 
elements of $\A$, it is not clear that this is a first order condition. Thus 
we need to show that such a theory $\ti{\Th}$ indeed exists, that it 
eliminates quantifiers, and that any model of $\Th_I$ can be embedded in a 
model of this kind.

The fact that the above condition is actually first order, follows from the 
following theorem of \Cite{bounds} (by the \noexpand\emph{degree} of a 
polynomial\index{degree of a polynomial} we mean the \emph{total} degree):

\begin{fact} \Label{fct}{bounds}
  Let $A$ be the polynomial algebra in $n$ variables over an arbitrary field, 
  $d$ a fixed degree. There is a degree $e$ depending only on $n$ and $d$ 
  (and not on the base field), such that for any $p_1,\dots,p_m\in{}A$ of 
  degree at most $d$:
  \begin{enumerate}
    \item
      For any $f\in{}A$ of degree at most $d$, if $f$ is in the ideal 
      generated by the $p_i$, then $f=\sum{}h_ip_i$ for $h_i$ of degree at 
      most $e$.
      
    \item
      The module of tuples $(s_1,\dots,s_m)$ such that $\sum{}s_ip_i=0$ is 
      generated by tuples of elements of degree at most $e$.

  \end{enumerate}

  More generally, the same results hold when $A$ is replaced by $A^k$.  Here, 
  the degree of $(t_1,\dots,t_k)\in{}A^k$ is the maximum of the degrees, and 
  $e$ depends also on $k$.
\end{fact}

\begin{remark}\Label{rmk}{algstr}
  \mbox{}
\begin{enumerate}
  \item
    For the polynomial algebra, the set of polynomials of a given degree 
    forms, in a natural way, a definable set. For the more general algebra 
    $\A$ we may define the degree of an element $r$ to be the minimal degree 
    of a pre-image of $r$ in $K[T_1,\dots,T_n]$. A priori, it is not clear 
    that the set of elements of a given degree in $\A$ is again definable, 
    since an element is represented in more than one way as a polynomial.  
    However, since (according to fact~\ref{fct:bounds}) membership in $I$ is 
    a first order property (of the coefficients), we have formulas whose free 
    variables represent an element of $\A$ of a given degree. Alternatively, 
    for the purpose of describing members of $\A$ we may assume that $\A$ is 
    actually the polynomial algebra, since $I$ only appears as a condition on 
    the modules.

  \item
    Elements of $K^m$ will usually be considered as coefficients of 
    polynomials in the $T_i$. This means that we fix an order on the 
    monomials in the $T_i$, and for $\bar{x}\in K^m$, $x_i$ is the 
    coefficient of the $i$-th monomial. Multiplication of polynomials induces 
    an operation $*:K^m\times K^l\MapsTo K^N$.

    The same is true for elements of the $G^i$: if two such elements $p,q$ 
    corresponds to vector spaces $V_p$ and $V_q$ of polynomials in the $T_i$, 
    $p*q$ corresponds to the image of $V_p\tensor{}V_q$ in the polynomial 
    algebra.

    Sometimes, instead of thinking of a tuple as a polynomial, we think of it 
    as a tuple of polynomials (it will be clear from the context.) In that 
    case, multiplication (by a polynomial or one vector space) is done 
    term-wise.

  \item
    Here is an instance of the above notation: Let $J$ be an ideal in $\A$. 
    Then $J$ is finitely generated; let $p$ be the vector space generated by 
    a finite set of generators. It follows from fact~\ref{fct:bounds}, that 
    given a degree $d$ there is a degree $e$ such that, setting $q=K^e$, the 
    set of elements of degree $d$ in $J$ is precisely the set of elements of 
    degree $d$ in $q*p$.

  \item
    Some more notation: $D_m(v,\bar{v})$ will denote 
    $D(T^{\bar{i}}v,\dots,v,\bar{v})$, where the $\dots$ stands for all 
    monomials of total degree at most $m$ (with the prescribed order).

  \item
    Recall from the case $n=0$ that over a good substructure, the type of 
    $D(v,\bar{v})$ determines the type of $D(v,\bar{u})$ whenever both are not 
    $0$ and $D(\bar{v})=D(\bar{u})=0$. The passage to a good substructure is 
    done precisely as in the previous case.
\end{enumerate}
\end{remark}

The fact that the condition on the $v_i$ is first order now follows from the 
second item of fact~\ref{fct:bounds}, since it is enough to state the 
conditions on the $f_i$ for generators of the tuples $(t_i)$. The fact that 
the condition on the $g_i$ is first order follows from the first item of 
fact~\ref{fct:bounds}.

Using the last point of remark~\ref{rmk:algstr}, we may obtain a description 
of the types:

\begin{claim}[Description of types, general case]\Label{clm}{types}
  For any quantifier free $1$-type $\p(v)$, either there are $m$ and 
  $\bar{v}$ such that $\p$ is determined by the formulas in it of the form 
  $\varphi(D_m(v,\bar{v}))$ (where $\varphi$ does not involve any module 
  stuff), or $\p$ is the unique quantifier free type determined by the set of 
  formulas $D_m(v,\bar{v})=0$ for all $m$ and $\bar{v}$.
\end{claim}

\begin{proof}

  Let $N$ be a model realising $\p$, $v\in{}M(N)$ a realisation. Since we are 
  working over a good substructure $N_0$, we may view 
  $K(N)\otimes_{K(N_0)}M(N_0)$ as a sub $\A$-module of $M(N)$. Let $J$ be the 
  ideal in $\A$ of elements $f$ such that 
  $fv\in{}K(N)\otimes_{K(N_0)}M(N_0)$. If this ideal is $0$, we are in the 
  second case. Otherwise, let $f_1,\dots,f_n$ generate $J$, and let 
  $\bar{v}_i\in{}M(N_0)$, for $i$ between $1$ and $n$, be bases for the 
  minimal $K(N_0)$ subspace containing $f_iv$.  We set $m$ to be the maximum 
  of the degrees of the $f_i$ and $\bar{v}=(\bar{v}_1,\dots,\bar{v}_n)$. A 
  different choice of $N$ and $v$ will result choosing $f_i$ of the same 
  form, with coefficients satisfying the same type over $K(N_0)$. Thus $m$ 
  and $\bar{v}$ do not depend on the choice of $N$ and $v$.
  
  Let $V_i$ be the vector space spanned by $\bar{v_i}$. Let $g\in\A$ be such 
  that $\p$ says that $D(gv,u_1,\dots,u_k)\ne{}0$ for some module elements 
  $u_i$. Then $g=\sum{}h_if_i$ for some $h_i\in{}\A$. Since the base is a 
  structure, applying the operators $T_i$ to elements of $V_i$ is well 
  defined. If $\psi(v)=\varphi(D(gv,\bar{u}))$ is a formulas in $\p$, 
  consider the definable set
  \begin{equation*}
    \{(w_1,\dots,w_n)\in V_1\oplus\dots\oplus V_n | \varphi(D(\sum h_iw_i, 
    u_1,\dots,u_k))\}
  \end{equation*}
  (The spaces $V_i$ are represented by tuples of field elements, so this is a 
  subset of the $G$ sorts.) Since we are over a good substructure, the type 
  of this space is determined by the base.  Also, $v$ satisfies $\psi$ if and 
  only if $(f_1v,\dots,f_nv)$ belongs to this set. But this is determined by 
  the type of $D_m(v,\bar{v})$ (over the base $N_0$)
\end{proof}

As in the case $n=0$, the result we seek easily follows from this:
\begin{thm} \Label{thm}{eq}
  Let $\tilde{\Th}$ be the theory extending $\Th$, and stating, in addition, 
  that for any $f_1,\dots,f_m,g_1,\dots,g_k\in\A$ and any 
  $v_1,\dots,v_m,\bar{u}\in{}M$ such that, for any $i$, $g_i$ is not in the 
  ideal generated by $f_1,\dots,f_m$, and such that for any 
  $t_1,\dots,t_m\in\A$, $\sum{}t_if_i=0$ implies $\sum{}t_iv_i=0$, the 
  formula
  \begin{equation*}
    \bigand{i} f_ix=v_i\land\bigand{i} D(g_ix,\bar{u})=0
  \end{equation*}
  has a solution $x$.

  Then $\tilde{\Th}$ eliminates quantifiers.
\end{thm}

\begin{proof}
  Using criterion~\ref{crt:eq}, and the above claim, we need to show that 
  given a good substructure $M_0$, and a quantifier free type $\p$ over 
  $K(M_0)$ and $G(M_0)$, we may satisfy $\p(D_m(v,\bar{v}))$ in any theory 
  extending $\Th_{M_0}$.
  
  Since $\p$ is a type in the $G$ sorts over $K(M_0)$, it follows from 
  section~\ref{sec:eq0} that $\p$ is consistent. Let $p$ satisfy $\p$. Again, 
  by the case $n=0$, we may assume that $p$ is in $M_0$, and we may extend 
  the field so that $p$ corresponds to some subspace of $K^l$. This means 
  that satisfying $\p(D_m(v,\bar{v}))$ amounts to satisfying conditions of 
  the form
  \begin{gather*}
    fv=\sum x_iv_i\\
    gv\notin \Gen{v_j}\\
  \end{gather*}
  Where $f$ is an element of $\A$ and $x_i\in K$. Since $\p$ was consistent 
  to start with, the conditions appearing in the axioms are satisfied for any 
  set of conditions like this that appears in $\p$. Hence, the axioms imply 
  that these equations have a solution.
\end{proof}

\begin{cor}
  $\tilde{\Th}$ is $\w$-stable
\end{cor}

\begin{proof}
  from the theorem, by counting the types.
\end{proof}
\begin{remark}\label{rmk:noncom}\index{algebra!non-commutative}
  Following through the proofs, one sees that they work just as well with the 
  commutativity assumption on the generators replaced by some other axioms, 
  provided that the resulting algebra is (left) Noetherian, and the class of 
  modules satisfying the solvability conditions is first order. In 
  particular, using the more general version of fact~\ref{fct:bounds}, we see 
  that the same result holds for algebras finite over their centre (where the 
  field is contained in the centre.)
\end{remark}

The last thing is to prove that any module over $\A$ (considered as a model 
of $\Th$) can be embedded in a model of $\ti\Th$. First note that the axioms 
can be split into two parts:
\begin{enumerate}
  \item\label{itm:injective}
    There is a solution for any finite set of equations $f_iv=u_i$, provided 
    that if $\sum t_if_i=0$ then $\sum t_iu_i=0$.

  \item
    There is a solution for any finite set of formulas $f_jv=0$, 
    $g_iv\notin{}U_i$ (where $U_i$ is a finite dimensional vectors space), 
    provided that no $g_i$ is in the ideal generated by the $f_j$.
\end{enumerate}

This is true since a solution of a general set of equations of the type 
considered is the sum of a solution of the corresponding equations of the 
first kind, and of the second kind.

We claim:
\begin{claim}
  Let $\A$ be a Noetherian ring. An $\A$ module $M$ satisfies condition 
  (\ref{itm:injective}) above (for $f_i,t_i\in\A$ and $v,u_i\in{}M$) if and 
  only if $M$ is injective.\index{module!injective}
\end{claim}

\begin{proof}
  Let $M$ be an injective $\A$-module, let $U\subseteq{}M$ be the submodule 
  generated by the $u_i$, and let $V=(U\oplus\A)/<f_i-u_i>$. The condition
  \begin{equation*}
    \sum t_if_i=0 \implies \sum t_iu_i=0
  \end{equation*}
  is equivalent to the map $U\MapsTo V$ being injective. Therefore, the 
  inclusion map of $U$ in $M$ extends to $V$, and the image of $1\in V$ in 
  $M$ is a solution.

  Conversely, by a result of Baer (cf.~\Cite{eisenbud}), it is enough to 
  check the condition of injectivity for the inclusion of an ideal $I$ in 
  $\A$.  Since $\A$ is Noetherian, $I$ is finitely generated, say by $f_i$.  
  Let $u_i$ be the images of $f_i$ in $M$. Then $f_i$, $u_i$ satisfy the 
  assumption of (\ref{itm:injective}), so there is some $v\in M$ such that 
  $f_iv=u_i$ for all $i$. Now, the map from $\A$ to $M$ that takes $1$ to $v$ 
  extends the given map.
\end{proof}

Regarding the other condition, consider the module $M=\prod{}M_I$, where 
$M_I=(\A/I)^{\w}$, for $I$ an arbitrary ideal of $\A$ (so the product is over 
all ideals.) We claim that any module containing $M$ satisfies the second 
condition. To see this, we first note that it is enough to show that $M$ 
itself satisfies the condition. Indeed, given arbitrary finite dimensional 
vector spaces $U_i$ in a module containing $M$, any solution in $M$ to the 
problem, with $U_i$ replaced by $U_i\cap{}M$ will solve the original problem.

For $M$ itself, let $I$ be the ideal generated by the $f_j$. For the same 
reason as before, it is enough to find a solution in $M_I$ (note that the 
condition is non-trivial only if $I$ is a proper ideal.)

Now, any element of $M_I$ is a solution to the equations. Thus we only need 
to satisfy the inequalities. Since the $g_i$ are not in $I$, they are 
non-zero in each $\A/I$. Hence, almost all of the unit vectors in $M_I$ 
satisfy the inequalities. This solves the problem.

We now can prove:

\begin{claim}\Label{clm}{embed}
  Any module over $\A$ embeds into a model of $\ti\Th$.
\end{claim}

\begin{proof}
  Let $N$ be any module. Then $N\oplus{}M$ can be embedded into some 
  injective module $I$ (where $M$ is the module constructed above.) Then $I$ 
  contains $N$, satisfies the first condition since it is injective, and 
  satisfies the second condition since it contains $M$.
\end{proof}

Finally, combining theorem~\ref{thm:eq} and claim~\ref{clm:embed}, we get:

\begin{cor}
  The theory $\ti\Th$ is the model completion of the theory $\Th=\Th_I$.
\end{cor}

Most definability results coming from algebra (such as fibre dimension) are 
concerned with \emph{finitely generated} modules. The following example shows 
the theory of a finitely generated module is far from having quantifier 
elimination. Therefore, such definability results can not be derived directly 
by considering the theory of the module, but should probably be obtained by 
interpreting the module in our theory $\ti\Th$.

\begin{example}
  Let $\A=K[T]$, the polynomial algebra in one variable over a field $K$ of 
  characteristic $0$, and consider $M=\A$ as a module over itself. We will 
  show that the semi-ring of natural numbers can be interpreted in this 
  theory.

  For elements $v$ of $M$ and $x$ of $K$, denote by $v(x)=0$ the formula 
  $\E{}u((T-x)u=v)$. Now consider the formula:

  \begin{equation*}
    \begin{split}
      \phi(v,y)&=\\
      & v\neq 0\land v(0)=0\land \\
      & \forall x(v(x)=0\implies (x=y\lor v(x+1)=0))
    \end{split}
  \end{equation*}

  We claim that for any $y$, the fibre $\phi(v,y)$ is non-empty (in $M$) if 
  and only if $y$ is a natural number. Indeed, assume that $y$ is not in 
  $\Z$, and that $v$ satisfies $\phi(v,y)$. Then $v$ is a non-zero element 
  that is divisible by $T+n$ for any natural $n$. There is no such element in 
  $K[T]$ (here we use that the characteristic is $0$.) Conversely, when $y$ 
  is natural, the element $T(T-1)\dots(T-y)$ satisfies the formula.

  The conclusion is that $\E{}v\in{}M(\phi(v,y))$ defines the set of natural 
  numbers. The ring operations are automatically defined, since this formula 
  actually defines the copy of the natural numbers contained in $K$.
\end{example}

This example holds more generally: If $\A$ is any finitely generated algebra 
over a field $K$ of characteristic $0$, and $M$ is any finitely generated 
module over $\A$ of infinite dimension over $K$, there is dominant map of 
$spec(\A)$ to the affine line, that makes $M$ into a $K[T]$ module, which, 
after a localisation becomes free (and in particular, torsion free.) The fact 
that $M$ is infinite dimensional means that the support of $M$ has dimension 
at least $1$, so that this map can be chosen so that the resulting free 
module is non-zero. Now we may repeat the above example to interpret (all but 
finitely many of) the natural numbers.


\bibliographystyle{amsplain}

\bibliography{../../bibtex/mr}

\end{document}